\documentclass[12pt]{article}
\usepackage{amsfonts}
\usepackage{amsmath}
\usepackage{latexsym,amssymb}
\usepackage[T1]{fontenc}
\usepackage[english]{babel}
\usepackage{amssymb}
\usepackage{theorem}
\usepackage{eufrak}

\newtheorem{theorem}{Th\'eor\`eme}[section]
\newtheorem{proposition}[theorem]{Proposition}
\newtheorem{corollary}[theorem]{Corollaire}
\newtheorem{e-definition}[theorem]{Definition\rm}

\newtheorem{lemma}[theorem]{Lemme}

\begin{document}

\author{Mohamed Taba\^a 
\quad \\
{\it \small
Faculté des Sciences, Département de Mathématiques,
B.P. 1014 Rabat Maroc
 }}

\title{Sur le produit tensoriel d'alg\`{e}bres\footnote{Le théorème 2.1, le Corollaire 2.2, le Corollaire 2.4 et la Proposition 2.6 de la première version de cet article (Arxiv: 1304.5395v1 [math.AC]) figurent dans le projet de note intitulé "Sur le produit tensoriel d'algèbres" que l'auteur avait soumis aux CRAS de Paris (n$^o$ CRMATHEMATIQUE-D-12-00204) le 07 Juin 2012.}}
\date{}
\maketitle
\begin{abstract}
Let $\sigma:A\rightarrow B$ and $\rho:A\rightarrow C$\ be
two homomorphisms of noetherian rings such that $B\otimes_{A}C$ is a
noetherian ring. we show that if $\sigma$ is a regular (resp. complete intersection, resp.
Gorenstein, resp. Cohen-Macaulay, resp. $%
(S_{n}$), resp. almost Cohen-Macaulay) homomorphism,  so is $\sigma\otimes I_{C}$ and the converse is true if $\rho$ is faithfully flat.
 We deduce the transfert of the previous properties of $B$ and $C$ for $B\otimes_{A}C$, and then for the completed tensor product $B\hat{\otimes}_{A}C$. If  $B\otimes_{A}B$ is noetherian and $\sigma$ is flat, we give a necessary and sufficient condition to
 $B\otimes_{A}B$ be a regular ring.
 \end{abstract}

\section{Introduction}

Tous les anneaux consid\'{e}r\'{e}s sont suppos\'{e}s commutatifs et
unitaires. Les notations sont celles de \cite{Grot1}.

Rappelons (\cite{Grot1},7.3.1) que si $\sigma :A\rightarrow B$ un homomorphisme d'anneaux noeth%
\'{e}riens, on dit que $\sigma $ est r\'{e}gulier s'il est plat et si pour
tout id\'{e}al premier $\mathfrak{p}$\ de $A$ la $k(\mathfrak{p})$-alg\`{e}%
bre $B\otimes _{A}k(\mathfrak{p})$ est g\'{e}om\'{e}triquement r\'{e}guli%
\`{e}re, et que $\sigma $ est d'intersection compl\`{e}te (resp.  de
Gorenstein, resp. de Cohen-Macaulay, resp.($S_{n}$), resp. presque de Cohen-Macaulay) s'il
est plat et si pour tout id\'{e}al premier $\mathfrak{p}$\ de $A$ l'anneau $%
B\otimes _{A}k(\mathfrak{p})$ est d'intersection compl\`{e}te (resp. de
Gorenstein, resp. de Cohen-Macaulay, resp. v\'{e}rifie ($%
S_{n}$), resp. presque de Cohen-Macaulay).

Dans ce qui suit nous montrons, que si $\sigma :A\rightarrow B$ et $\rho
:A\rightarrow C$\ sont deux homomorphismes d'anneaux noeth\'{e}riens tels
que $B\otimes _{A}C$ soit un anneau noeth\'{e}rien, alors  $\sigma\otimes I_{C}$  est r\'{e}gulier (resp.  d'intersection compl\`{e}te, resp. de
Gorenstein, resp. de Cohen-Macaulay, resp. ($S_{n}$), resp. presque de Cohen-Macaulay) si $\sigma $ l'est; et
que la r\'{e}ciproque est vraie si $\rho $ est fid\`{e}lement plat. On en d\'{e}duit, en particulier, que si
$\sigma$ est plat, alors $B\otimes_{A}C$ est un anneau d'intersection compl\`{e}te (resp.  de
Gorenstein, resp. de Cohen-Macaulay) si $B$ et $C$ le sont, et il est presque de Cohen-Macaulay si l'un des anneaux, $B$ ou $C$, est de Cohen-Macaulay et l'autre est presque de Cohen-Macaulay.\

Si $A$ est un corps,  on retrouve le Th\'{e}or\`{e}me 2 de \cite{Watanabe} si $B$
et $C$ sont des anneaux de Gorenstein, le Th\'{e}or\`{e}me 2.1 de  \cite{abou1}
si $B$ et $C$ sont des anneaux de Cohen-Macaulay et le Th\'{e}or\`{e}me 6 de
\cite{Tousi}  si $B$ et $C$ sont des anneaux d'intersection compl\`{e}te (resp. v%
\'{e}rifient ($S_{n}$)).

Comme application nous montrons que si $\sigma$ et $\rho$ sont deux
homomorphismes locaux d'anneaux locaux noeth\'{e}riens, et si le corps r\'{e}%
siduel de $C$ est de rang fini sur celui de $A$, alors le produit tensoriel
compl\'{e}t\'{e} $B\hat{\otimes}_{A}C$ est r\'{e}gulier\ si l'homomorphisme $%
\sigma$ est formellement lisse et $C$ est r\'{e}gulier, et  il est d'intersection compl\`{e}te (resp. de Gorenstein,
resp. de Cohen-Macaulay) si $\sigma$
est plat et $B$ et $C$ le sont, et il est presque de Cohen-Macaulay si $\sigma$
est plat et si l'un des anneaux, $B$ ou $C$, est de Cohen-Macaulay et l'autre est presque de Cohen-Macaulay.

Si $\sigma$ est  plat et $B\otimes_{A}B$ est un anneau noeth\'{e}rien, nous montrons que $B\otimes_{A}B$ est r\'{e}gulier si et seulement
si $B$ est r\'{e}gulier et $\sigma$\ est r\'{e}gulier.

Dans toute la suite nous utilisons librement les r\'{e}sultats de \cite{Mat} et de \cite{Avramov1},
 et l'homologie d'Andr\'{e}-Quillen telle qu'elle est d\'{e}finie dans  \cite{andre}.

\section{R\'{e}sultats}

\begin{proposition}
Soit $\sigma:A\rightarrow B$ un homomorphisme d'anneaux noeth\'{e}%
riens. Les propri\'{e}t\'{e}s suivantes sont \'{e}quivalentes:

\begin{enumerate}
\item[i)] L'homomorphisme $\sigma$ est r\'{e}gulier (resp. d'intersection compl\`{e}te).

\item[ii)] L'homomorphisme $\sigma$ est plat et $H_{1}(A,B,k(\mathfrak{q}))=0$ (resp. $H_{2}(A,B,k(\mathfrak{q}%
))$ $=0$) pour tout id\'{e}al premier $\mathfrak{q}$\ de $B$.
\end{enumerate}
\end{proposition}

\noindent\textbf{D\'{e}monstration.} Cas r\'{e}gulier: cf. \cite{andre} Suppl\'ement Th\'{e}or\`{e}me 30.

Cas d'intersection compl\`{e}te: (cf. \cite{Marot} ) On utilise \cite{andre}. Soit $%
\mathfrak{q}$\ id\'{e}al premier\ de $B$ et $\mathfrak{p}=\sigma ^{-1}(%
\mathfrak{q})$. D'apr\`{e}s le Corollaire 5.27, la Propositon 4.54, la suite
exacte associ\'{e}e aux homomorphismes $k(\mathfrak{p})\rightarrow B_{%
\mathfrak{q}}/\mathfrak{p}B_{\mathfrak{q}}\rightarrow \ k(\mathfrak{q})$ et
d'apr\`{e}s la Proposition 7.4, on a $H_{2}(A,B,k(\mathfrak{q})\cong
H_{3}(B_{\mathfrak{q}}/\mathfrak{p}B_{\mathfrak{q}},k(\mathfrak{q}),k(%
\mathfrak{q}))$; l'\'{e}quivalence r\'{e}sulte
donc de la Proposition 6.27.\
\begin{lemma}
Soient $\sigma:A\rightarrow B$ et $\rho:A\rightarrow C$\ deux homomorphismes
d'anneaux, $\mathfrak{Q}$ un id\'{e}al premier de $B\otimes _{A}C$ et $\mathfrak{q}=(I_{B}\otimes\rho)^{-1}(\mathfrak{Q})$. Si $\sigma$ est plat, alors on l'isomorphisme suivant:
\[
 H_{n}(A,B,k(\mathfrak{q}))\otimes_{k(%
\mathfrak{q})}k(\mathfrak{Q})\cong H_{n}(C,B\otimes_{A}C,k(\mathfrak{Q}))
\]
\end{lemma}
\noindent\textbf{D\'{e}monstration.} En effet, d'apr\`es le Lemme 3.20 de \cite{andre} on a $H_{n}(A,B,k(\mathfrak{q}))\otimes_{k(%
\mathfrak{q})}k(\mathfrak{Q})\cong H_{n}(A,B,k(\mathfrak{Q}))$ et d'apr\`es la Proposition 4.54 de \cite{andre} on a $H_{n}(A,B,k(\mathfrak{Q}))$ $ \cong H_{n}(C,B\otimes_{A}C,k(\mathfrak{Q}))$; d'o\`u le Lemme.

 \begin{theorem}
 Soient $\sigma:A\rightarrow B$ et $\rho:A\rightarrow C$\ deux homomorphismes
d'anneaux noeth\'{e}riens. \ On suppose que $B\otimes_{A}C$ est un anneau
noeth\'{e}rien. Alors:

a) Si $\sigma$ est r\'{e}gulier, il en est de m\^{e}me de $\sigma\otimes I_{C}:C\rightarrow$ $B\otimes_{A}C$; la r\'{e}ciproque est vraie si $\sigma$ est plat et $^a\rho$ est surjective.

b) Si les fibres de $\sigma$ sont des anneaux d'intersection compl\`ete (resp.  de Gorenstein, resp. de Cohen-Macaulay, resp. v\'erifient $(S_{n})$), il
en est de m\^{e}me de celles de $\sigma\otimes I_{C}$;
la r\'{e}ciproque est vraie si $^a\rho$ est surjective.
\end{theorem}

\noindent\textbf{D\'{e}monstration.} a) Supposons que $\sigma$ est un homomorphisme r\'egulier , alors il est plat et par suite $\sigma\otimes I_{C}$ est plat. L'implication r\'esulte alors de la Proposition pr\'ec\'edente en tenant compte du Lemme.

R\'eciproquement, d'apr\`es Proposition (\textbf{I},3.6.1) de \cite{Grot2}, l'application  $^a(\sigma\otimes I_{C})$ est surjective. La r\'eciproque r\'esulte aussi de la Proposition pr\'ec\'edente en tenant compte du Lemme  .\

b) i) Supposons d'abord que $\sigma$ est un homomorphisme d'intersection compl\`ete, le m\^eme raisonnement que dans le cas pr\'ec\'edent montre que
$\sigma\otimes I_{C}$ est un homomorphisme d'intersection compl\`ete.

ii) Supposons maintenant que les fibres de $\sigma$ sont des anneaux d'intersection compl\`{e}te  (resp.  de Gorenstein, resp. de
Cohen-Macaulay, resp. v\'erifient $(S_{n})$ ). Posons $D=B\otimes_{A}C$ et\ soit $%
\mathfrak{r}$ un id\'{e}al premier de $C$. L'anneau $D\otimes_{c}k(\mathfrak{%
r})=$ $(B\otimes_{A}C)\otimes_{C}k(\mathfrak{r})$ est isomorphe \`{a} $%
B\otimes_{A}k(\mathfrak{r})$. Soit $\mathfrak{p}=\rho ^{-1}(\mathfrak{r})$.
Donc $D\otimes_{C}k(\mathfrak{r})$ est isomorphe \`{a} $(B\otimes _{A}k(%
\mathfrak{p}))\otimes_{k(\mathfrak{p})}k(\mathfrak{r})$. Comme
l'homomorphisme $k(\mathfrak{p})\rightarrow k(\mathfrak{r})$ est
d'intersection compl\`{e}te, il r\'{e}sulte du cas pr\'{e}c\'{e}dent appliqu\'{e}
aux homomorphismes $k(\mathfrak{p})\rightarrow k(\mathfrak{r})$ et $k(%
\mathfrak{p})\rightarrow B\otimes_{A}k(\mathfrak{p})$ que l'homomorphisme $%
B\otimes_{A}k(\mathfrak{p})\rightarrow (B\otimes_{A}k(\mathfrak{p}%
))\otimes_{k(\mathfrak{p})}k(\mathfrak{r})$ est d'intersection compl\`{e}te.
On en d\'{e}duit que l'homomorphisme $B\otimes_{A}k(\mathfrak{p})\rightarrow
D\otimes_{C}k(\mathfrak{r})$ est d'intersection compl\`{e}te (resp.  de Gorenstein, resp. de Cohen-Macaulay,
resp. $(S_{n})$) et que par suite $D\otimes_{C}k(\mathfrak{r})$ est un
anneau d'intersection compl\`{e}te (resp.  de Gorenstein, resp. de Cohen-Macaulay, resp. v\'{e}rifie  $(S_{n})$ ).

R\'eciproquement, soit $\mathfrak{p}$ un id\'{e}al premier de $A$ et soit  $\mathfrak{r}$ un id\'eal premier de $C$ tel que $\mathfrak{p}=\rho^{-1}(\mathfrak{r})$.
L'homomorphisme $k(\mathfrak{p})\rightarrow k(\mathfrak{r})$ est fid\`element plat, il en est de m\^eme de l'homomorphisme  $B\otimes_{A}k(\mathfrak{p})\rightarrow
D\otimes_{C}k(\mathfrak{r})$. Donc  $B\otimes_{A}k(\mathfrak{p})$
est un anneau d'intersection compl\`{e}te (resp.  Gorenstein,  resp. de Cohen-Macaulay, resp. v\'{e}%
rifie  $(S_{n})$ ). \\

\noindent {\it Remarques.} i) Si l'homomorphisme $\rho$ est fid\`{e}lement plat alors l'application $^a\rho$ est surjective et si de plus $\sigma\otimes I_{C}$ est plat alors $\sigma$ est plat.\

ii) Dans \cite{Avramov2}, \cite{Avramov3}, \cite{Avramov4}, on trouve des r\'{e}sultats sur le changement de base pour les homomorphismes qu'ils ont d\'{e}fini.

\begin{proposition}
Soient $\sigma:A\rightarrow B$ et $\rho:A\rightarrow C$\ deux homomorphismes
d'anneaux noeth\'{e}riens. \ On suppose que $B\otimes_{A}C$ est un anneau
noeth\'{e}rien et que $\sigma$ est r\'{e}gulier (resp.  d'intersection compl%
\`{e}te, resp. de Gorenstein, resp. de Cohen-Macaulay, resp. $\ (S_{n})$).
Si $C$ est un anneau r\'{e}gulier (resp.  d'intersection compl\`{e}te, resp.
de Gorenstein, resp. de Cohen-Macaulay, resp. v\'{e}rifie $(S_{n})$ ) il en
 est de m\^{e}me de $B\otimes_{A}C$; la r\'eciproque est vraie si $\sigma$ est fid\`element plat.
\end{proposition}
\noindent\textbf{D\'{e}monstration.} D'apr\`es le Th\'eor\`eme pr\'ec\'edent, l'homomorphisme $\sigma\otimes I_{C}$ est r\'egulier (resp.  d'intersection compl%
\`{e}te, resp. de Gorenstein, resp. de Cohen-Macaulay, resp. $\ (S_{n})$); d'o\`u l'implication. La r\'eciproque r\'esulte du fait que $\sigma\otimes I_{C}$ est fid\`element plat.\\

On en d\'{e}duit que si $k$ est un corps et si $B\otimes_{k}C$ est un anneau
noeth\'{e}rien alors $B\otimes_{k}C$ v\'{e}rifie  $%
(S_{n}) $ si B et C la v\'{e}rifient.

\begin{corollary}
Soit $k$ un corps. On suppose que $B\otimes _{k}C$ est un anneau noeth\'{e}rien et
que pour tout id\'{e}al maximal $\mathfrak{n}$ de $C$, $k(\mathfrak{n})$ est
s\'{e}parable sur $k$. Si $B$ et $C$ sont r\'{e}guliers alors $B\otimes
_{k}C$ est r\'{e}gulier.
\end{corollary}

\noindent\textbf{D\'{e}monstration.} Pout tout id\'{e}al maximal $\mathfrak{n}$ de $C,\ \ k(%
\mathfrak{n})$ est  s\'{e}parable sur $k$ et $C_{\mathfrak{n}}$ est r\'{e}%
gulier, donc $C_{\mathfrak{n}}$ est g\'{e}om\'{e}triquement r\'{e}guli\`{e}%
re sur $k$.  Le r\'{e}sultat d\'{e}coule donc de la Proposition pr\'ec\'edente puisque l'
homomorphisme $k\rightarrow C$ est r\'{e}gulier.\\

Si $C$ est r\'egulier, il peut se faire que $k(\mathfrak{n})$ soit
s\'{e}parable sur $k$ pour tout id\'{e}al maximal $\mathfrak{n}$
de $C$, sans que $k(\mathfrak{r})$ soit  s\'{e}parable sur $k$ pour tout
id\'{e}al premier $\mathfrak{r}$ de $C$ [6, Theorem 2.11].
En effet, soient $k$ un corps non parfait de caract\'eristique $p>0$, $a$ un \'el\'ement de $k-k^p$,
 $A$  l'anneau de polyn\^omes $k[X,Y]$ et
$C$ l'anneau local de $A$ en l'id\'{e}al maximal engendr\'{e} par $X$ et $Y$. L'anneau $C$ est r\'{e}gulier
et son corps r\'{e}siduel $k(X,Y)$ est s\'{e}parable sur $k$. Soit $f$ le polyn\^ome $Y^p-aX^p$.
Comme $a\notin k$, $f$ est irr\'{e}ductible dans $A$. Notons $\mathfrak{p}$  l'id\'{e}al premier
de $A$ engendr\'{e} par $f$ et $\mathfrak{r}$ l'id\'{e}al premier  $\mathfrak{p}C$ de $C$. Montrons que
$ k(\mathfrak{r})$ n'est pas s\'{e}parable sur $k$. Le corps $ k(\mathfrak{r})$ s'identifie
canoniquement \`a $k(\mathfrak{p})$ donc il suffit de montrer que $k(\mathfrak{p})$ n'est pas
 s\'{e}parable sur $k$. Si $x$ et $y$ sont les images respectives de $X$ et $Y$ dans
 $A{/\mathfrak{p}}$, on a bien $(\frac{y}{x})^p \in k$ et $\frac{y}{x}\notin k$.  Donc
 $k(\mathfrak{p})$ n'est  pas s\'{e}parable sur $k$.\

\begin{proposition}
Soient $\ \sigma :A\rightarrow B$ et $\rho :A\rightarrow C$\ deux
homomorphismes d'anneaux noeth\'{e}riens. \ On suppose que $B\otimes _{A}C$
est un anneau noeth\'{e}rien. Si $\sigma $ est plat alors $B\otimes _{A}C$
est un anneau d'intersection compl\`{e}te (resp.  de Gorenstein, resp. de
Cohen-Macaulay) si B et C le sont.
\end{proposition}
\noindent\textbf{D\'{e}monstration.} Si $B$ est un anneau d'intersection compl\`{e}te (resp.  de Gorenstein, resp. de
Cohen-Macaulay) alors $\sigma$ est un homomorphisme d'intersection compl\`{e}te (resp.  de Gorenstein, resp. de
Cohen-Macaulay). Le r\'esultat d\'ecoule donc de la Proposition  2.4.

\begin{proposition}
Soient $\ \sigma :A\rightarrow B$ et $\rho :A\rightarrow C$\ deux
homomorphismes locaux d'anneaux locaux noeth\'{e}riens. On suppose que le
corps r\'{e}siduel $C/\mathfrak{n}$ de $C$ est de rang fini sur le corps r%
\'{e}siduel $A/\mathfrak{m}$ de A.\newline
a) Si l'homomorphisme $\sigma $ est formellement lisse et $C$ est r\'{e}%
gulier alors l'anneau semi-local $B\hat{\otimes}_{A}C$ est r\'{e}gulier.%
\newline
b) Si $\sigma $ est plat alors $B\hat{\otimes}_{A}C$ est un anneau
d'intersection compl\`{e}te (resp. de Gorenstein, resp. de Cohen-Macaulay)
si $B$ et $C$ le sont.
\end{proposition}

\noindent \textbf{D\'{e}monstration.} Cas o\`{u} $B$ est complet. On utilise
la Proposition (\textbf{0},7.7.10) de \cite{Grot2}: Posons $E=B\hat{\otimes}_{A}C$ . D'apr%
\`{e}s i) $E$ est semi-local noeth\'{e}rien. Montrons que c'est un anneau r%
\'{e}gulier (resp. d'intersection compl\`{e}te , resp. de Gorenstein, resp.
de Cohen-Macaulay). Soit $\mathfrak{Q}$ un id\'{e}al maximal de $E$. D'apr%
\`{e}s ii) $\mathfrak{Q}$ est au dessus de $\mathfrak{n}$. Pour montrer que $%
E_{\mathfrak{Q}}$ est r\'{e}gulier (resp. d'intersection compl\`{e}te ,
resp. de Gorenstein, resp. de Cohen-Macaulay) il suffit de montrer que $(E/%
\mathfrak{n}E)_{\mathfrak{Q}}$ est r\'{e}gulier (resp. d'intersection compl%
\`{e}te , resp. de Gorenstein, resp. de Cohen-Macaulay) puisque $C$ l'est et
d'apr\`{e}s iii) $E$ est un $C$-module plat. D'apr\`{e}s ii) $E/\mathfrak{n}E
$ est isomorphe \`{a} $B\otimes _{A}(C/\mathfrak{n}).$ Donc b) r\'{e}sulte
de la Proposition 2.6. D'autre part $B\otimes _{A}(C/\mathfrak{n})$ est
isomorphe \`{a} $(B/\mathfrak{m}B)\otimes _{A/\mathfrak{m}}(C/\mathfrak{n})$
donc a) r\'{e}sulte de la Proposition 2.4 puisque l'homomorphisme $A/\mathfrak{m}%
\longrightarrow B/\mathfrak{m}B$ est r\'{e}gulier.

Cas g\'{e}n\'{e}ral. L'anneau $B\hat{\otimes}_{A}C$ s'identifie \`{a} $\hat{B%
}\hat{\otimes}_{A}C$. Il suffit d'appliquer le cas pr\'{e}c\'{e}dent aux
homomorphismes $A\overset{\sigma }{\longrightarrow }B\longrightarrow \hat{B}$
et $\rho $. Dans a) l'homomorphisme $A\overset{\sigma }{\longrightarrow }%
B\longrightarrow \hat{B}$ est formellement lisse et dans b) $\hat{B}$ v\'{e}%
rifie la m\^{e}me propri\'{e}t\'{e} que $B$.\\

\noindent {\it Exemples.}
Les deux exemples suivants montrent que les r\'esultats pr\'ec\'edents tombent en d\'efaut
si les homomorphismes ne sont pas plats.

On prend pour $A$ un anneau de valuation discr\`ete complet, $\pi$  une uniformisante de $A$ et $k$ son corps r\'esiduel.

i) Si $B=k$ et $C=A[[X]]/({X^2}-\pi)$ o\`u  $X$ est ind\'{e}termin\'{e}e sur $A$, alors $B$ et $C$ sont des anneaux locaux r\'eguliers et la $k$-alg\`ebre $B\otimes_{A}k$ est g\'{e}om\'{e}triquement r\'{e}guli\`ere
mais les anneaux  $B\otimes _{A}C$ et $B\hat{\otimes}_{A}C$ qui sont isomorphes d'apr\`es la Proposition (\textbf{0},7.7.9) de \cite{Grot2}, ne sont pas r\'eguliers car $B\otimes _{A}C$ est isomorphe \`a  $k[[X]]/(X^2)$.

ii) Si $B=A[[X,Y]]/({X^2}-\pi,XY)$ o\`u  $X$ et $Y$ sont deux ind\'{e}termin\'{e}es sur $A$ et $C=k$, les anneaux $B$ et $C$ sont des anneaux locaux d'intersection compl\`ete, mais les anneaux  $B\otimes _{A}C$ et $B\hat{\otimes}_{A}C$ ne sont pas des anneaux de Cohen-Macaulay car  ils sont isomorphes \`a  $k[[X,Y]]/(X^2,XY)$. \

\begin{proposition}
Soit $\sigma:A\rightarrow B$ un homomorphisme d'anneaux noeth\'{e}riens.
On suppose que $B\otimes_{A}B$ est un anneau noeth\'{e}rien et que $\sigma$
 plat. Alors les propri\'{e}t\'{e}s suivantes sont \'{e}%
quivalentes:\newline
i) L'anneau $B$ est r\'{e}gulier et l'homomorphisme $\sigma$\ est r\'{e}%
gulier.\newline
ii) L'anneau $B\otimes_{A}B$ est r\'{e}gulier.
\end{proposition}

\noindent\textbf{D\'{e}monstration.} $i)\Rightarrow ii)$ Cela r\'{e}sulte du
Corollaire 2.2.

$ii)\Rightarrow i)$ Supposons $B\otimes_{A}B$\ r\'{e}gulier. L'homomorphisme $%
\sigma$\ est plat donc $\sigma\otimes I_{B}$\ est fid\`{e}%
lememt plat et par suite $B$ est r\'{e}gulier. Montrons que $\sigma$\ est r\'{e}gulier. Soit $\mathfrak{q}$ un
id\'{e}al premier de B. Comme $H_{1}(A,B,k(\mathfrak{q}))\cong H_{1}(B,B\otimes_{A}B,k(\mathfrak{q}))$ il suffit de montrer que $H_{1}(B,B\otimes_{A}B,k(\mathfrak{q}))=0$. On a la suite exacte
\[
H_{2}(B\otimes_{A}B,B,k(\mathfrak{q}))\rightarrow H_{1}(B,B\otimes_{A}B,k(%
\mathfrak{q}))\rightarrow H_{1}(B,B,k(\mathfrak{q}))
\]
associ\'{e}e \`{a} la factorisation $p\circ(\sigma\otimes_{A} I_{B})=I_{B}$,
o\`{u} $p:B\otimes_{A}B\rightarrow B$ est l'homomorphisme canonique d\'{e}%
fini par $p(b\otimes b^{\prime})=bb^{\prime}$. D'apr\`{e}s \cite[Suppl\'{e}%
ment, Proposition. 32]{andre} \ on a $H_{2}(B\otimes_{A}B,B,k(\mathfrak{q}))=0$. Donc $H_{1}(B,B\otimes_{A}B,k(\mathfrak{q}))=0$ car on a $H_{1}(B,B,k(%
\mathfrak{q}))=0$.

\begin{corollary}
Soient $K$ un corps et $L$ une extension de $K$. On suppose que l'anneau $%
L\otimes_{K}L$ est noeth\'{e}rien. Alors $L\otimes_{K}L$ est un anneau r\'{e}%
gulier si et seulement si L est s\'{e}parable sur K.
\end{corollary}\
\noindent\textbf{D\'{e}monstration.} En effet, l'homomorphisme $K\rightarrow L$ est r\'egulier si et seulement si $L$ est une extension s\'eparable de $K$.

\section{Cas Presque Cohen-Macaulay}
Suivant (\cite{Kang},1.5) on dira qu'un anneau noeth\'erien $A$ est presque de Cohen-Macaulay si $dim(A_{\mathfrak{p}})\leq{prof(A_{\mathfrak{p}}})+1$ pour tout id\'eal premier $\mathfrak{p}$ de $A$.\

 Il est clair que si $A$ est presque de Cohen-Macaulay alors, pour tout id\'eal premier $\mathfrak{p}$ de $A$, $A_{\mathfrak{p}}$ est presque de Cohen-Macaulay; et d'apr\`es (\cite{Kang}, 2.6), $A$ est presque de Cohen-Macaulay, si et seulement si, $dim(A_{\mathfrak{m}})\leq{prof(A_{\mathfrak{m}}})+1$ pour tout id\'eal maximal $\mathfrak{m}$ de $A$.\

Le r\'esultat suivant est une variante de la Proposition 2.2 de \cite{Ion} qui distingue les anneaux presque de Cohen-Macaulay des anneaux consid\'er\'es dans le paragraphe pr\'ec\'edent.
\begin{lemma}
Soient $\sigma:A\rightarrow B$ un homomorphisme local d'anneaux locaux noeth\'{e}riens et $\mathfrak{m}$ l'id\'eal maximal de $A$. On suppose que
$\sigma:A\rightarrow B$ est plat. Alors les propri\'et\'es suivantes sont \'equivalentes:\newline
a) L'anneau $B$ est presque de Cohen-Macaulay.\newline
b) L'un des anneaux, $A$ ou $B/\mathfrak{m}B$, est de Cohen-Macaulay et l'autre est presque de Cohen-Macaulay.
\end{lemma}
\noindent\textbf{D\'{e}monstration.} L'\'equivalence r\'esulte des \'egalit\'es
\[
dimB=dimA+dim{B/\mathfrak{m}B}\newline
\]
\[
profB=profA+prof{B/\mathfrak{m}B}
\]
et de la double in\'egalit\'e
\[
profA\leq{dimA}\leq{profA+1}.
\]

L'exemple suivant donne une r\'eponse \`a la question 2.3 pos\'ee par Ionescu dans \cite{Ion}.

{\it Exemple.}
Soient $k$ un corps, $R$ l'anneau $k[X,Y]/(X^2,XY)$ et $S$
l'anneau $k[X,Y,U,V]/(X^2,XY,U^2,UV)$.
L'homomorphisme canonique $R\rightarrow S$ est plat car
 L'homomorphisme $R\rightarrow R\otimes_{k}R$ est plat et $S$ s'identifie \`a
 $R\otimes_{k}R$. Soient $A$ l'anneau local de $R$ en l'id\'eal maximal
 $(x,y)$ et $B$ celui de $S$ en l'id\'eal maximal $(x,y,u,v)$. L'homomorphisme induit $A\rightarrow B$
 est local et plat. On v\'erifie que $xu$ annule l'id\'eal $(x,y,u,v)$ et que $xu\neq0$. Donc $profB=0$
 et par suite $profA=prof{B/\mathfrak{m}B}=0$. D'autre part on a  $dimA=dim{B/\mathfrak{m}B}=1$ et $dimB=2$. Donc les anneaux $A$ et ${B/\mathfrak{m}B}$ sont des anneaux presque de Cohen-Macaulay mais $B$ ne l'est pas.\\

Pour les homomorphismes presque de Cohen-Macaulay on a:
\begin{proposition}
Soient $\sigma:A\rightarrow B$ et $\rho:A\rightarrow C$\ deux homomorphismes
d'anneaux noeth\'{e}riens. \ On suppose que $B\otimes_{A}C$ est un anneau
noeth\'{e}rien. Si les fibres de $\sigma$ sont presque de Cohen-Macaulay  il
en est de m\^{e}me de celles de $\sigma\otimes I_{C}$; la r\'{e}ciproque est vraie si $^a\rho$ est surjective..
\end{proposition}
\noindent\textbf{D\'{e}monstration.} Soient $\mathfrak{r}$ un id\'{e}al premier de $C$ et $\mathfrak{p}=\rho ^{-1}(\mathfrak{r})$. L'homomorphisme $k(\mathfrak{p})\rightarrow k(\mathfrak{r})$
est de Cohen-Macaulay donc, d'apr\`es le Th\'eor\`eme, l'homomorphisme $B\otimes_{A}k(\mathfrak{p})\rightarrow
({B\otimes_{A}C})\otimes_{C}k(\mathfrak{r})$ l'est aussi; l'implication r\'esulte alors du Lemme pr\'ec\'edent. La r\'eciproque en r\'esulte aussi puisque
l'homomorphisme est fid\`element plat.
\begin{proposition}
Soient $\ \sigma :A\rightarrow B$ et $\rho :A\rightarrow C$\ deux
homomorphismes d'anneaux noeth\'{e}riens. \ On suppose que $B\otimes _{A}C$
est un anneau noeth\'{e}rien et que $\sigma $ est plat. Si pour tout id\'eal maximal $\mathfrak{N}$ de $B\otimes _{A}C$
l'un des anneaux, $B_{\mathfrak{q}}$
 ou $C_{\mathfrak{r}}$, o\`u $\mathfrak{q}$ et $\mathfrak{r}$  sont les images r\'eciproques respectives dans $B$ et $C$, est de Cohen-Macaulay et l'autre est presque de Cohen-Macaulay alors l'anneau $B\otimes _{A}C$ est presque de Cohen-Macaulay.
\end{proposition}\
\noindent\textbf{D\'{e}monstration.} L'anneau ${(B\otimes _{A}C)}_{\mathfrak{N}}$ s'identifie \`a un anneau de fractions de $B_{\mathfrak{q}}\otimes _{A_{\mathfrak{p}}}C_{\mathfrak{r}}$, o\`u $\mathfrak{p}$ est l'image r\'eciproque de $\mathfrak{N}$ dans $A$. Donc on se ram\`ene au cas o\`u  l'un des anneaux, $B$ ou $C$, est de Cohen-Macaulay et l'autre est presque de Cohen-Macaulay.

Supposons que  $B$ est un anneau de Cohen-Macaulay (resp.  presque de Cohen-Macaulay ); alors dans ce cas, la conclusion d\'ecoule du Lemme puisque $\sigma$ est un homomorphisme de Cohen-Macaulay (resp.  presque de Cohen-Macaulay )
 et par suite, d'apr\`es le Th\'eor\`eme (resp. la Proposition pr\'ec\'edente ), $\sigma\otimes I_{C}$ l'est.\\

\begin{proposition}
Soient $\ \sigma :A\rightarrow B$ et $\rho :A\rightarrow C$\ deux
homomorphismes locaux d'anneaux locaux noeth\'{e}riens. On suppose que le
corps r\'{e}siduel de $C$ est de rang fini sur celui de A. Si $\sigma $ est plat alors l'anneau $B\hat{\otimes}_{A}C$ est presque de Cohen-Macaulay
si l'un des anneaux, $B$ ou $C$, est de Cohen-Macaulay et l'autre est presque de Cohen-Macaulay.
\end{proposition}
\noindent\textbf{D\'{e}monstration.} On raisonne comme dans la Proposition 2.7: Si $B$ est un anneau de Cohen-Macaulay (resp.  presque de Cohen-Macaulay ) on utilise la Proposition 2.6(resp.  3.3) puis le Lemme.\\

L'exemple suivant montre que les deux Propositions pr\'ec\'edentes tombent en d\'efaut si les deux anneaux ( locaux ) sont presque de Cohen-Macaulay.\

{\it Exemple.} Soient $k$ un corps, $B$ l'anneau local de $k[X,Y]/(X^2,XY)$ en l'id\'eal maximal
 $(x,y)$, $C$ celui de $k[U,V]/(U^2,UV)$  en l'id\'eal maximal $(u,v)$. Les anneaux $B$ et $C$ sont presque de Cohen-Macaulay et $B\otimes_{k}C$ est noeth\'erien. Notons $D$ l'anneau $B\otimes_{A}C$, $\mathfrak{r}$ l'id\'eal maximal ${(x,y)\otimes_{k}B+C\otimes_{k}(u,v)}$ de $D$,  $E$ l'anneau $B\hat{\otimes}_{A}C$ et $T$ l'anneau local de $k[X,Y,U,V]/(X^2,XY,U^2,UV)$ en l'id\'eal maximal $(x,y,u,v)$.

 i) L'anneau $D_{\mathfrak{r}}$ est isomorphe \`a $T$ donc $D$ n'est pas un anneau presque de Cohen-Macaulay.

 ii) L'anneau $E$ est complet donc $\mathfrak{r}E$ est contenu dans son radical, et $\mathfrak{r}E$ est un id\'eal maximal de $E$ car $E/\mathfrak{r}E$ est isomorphe \`a $D/\mathfrak{r}$. Donc $E$ est un anneau local et
 $\mathfrak{r}E$ est son id\'eal maximal. L'homomorphisme $D\rightarrow E$ est plat. Donc l' homomorphisme induit $D_{\mathfrak{r}}\rightarrow E$
 est local et plat. D'apr\`es le Lemme pr\'ec\'edent,  $E$ n'est pas un anneau presque de Cohen-Macaulay
 car $D_{\mathfrak{r}}$ ne l'est pas.\\












\end{document}